\documentclass[12pt]{extarticle}
\usepackage{amsmath, amsthm, amssymb, hyperref, color}
\usepackage[shortlabels]{enumitem}
\usepackage{graphicx}
\usepackage[all]{xypic}
\usepackage{makecell}
\usepackage[final]{pdfpages}
\setboolean{@twoside}{false}
\usepackage{pdfpages}
\usepackage{caption}
\usepackage{subcaption}
\usepackage{scalefnt}
\usepackage{verbatim}
\oddsidemargin \evensidemargin
\newtheorem{theorem}{Theorem}
\numberwithin{theorem}{section}
\newtheorem{proposition}[theorem]{Proposition}
\newtheorem{lemma}[theorem]{Lemma}
\newtheorem{corollary}[theorem]{Corollary}
\newtheorem{definition}[theorem]{Definition}

\newtheorem{remark}[theorem]{Remark}
\newtheorem{example}[theorem]{Example}

\newtheorem{question}[theorem]{Question}

\newcommand{ \crk }{{\mathrm{rk}_{\mathbb{C}}}}
 \date{}

\title{\textbf{Real rank boundaries and loci of forms}}

\author{ Emanuele Ventura}

\begin{document}

\maketitle

\begin{abstract} 

\noindent In this article we study forbidden loci and typical ranks of forms with respect to the embeddings of $\mathbb P^1\times \mathbb P^1$ given by the line bundles $(2,2d)$. We introduce the Ranestad-Schreyer locus corresponding to supports of non-reduced apolar schemes. We show that, in those cases, this is contained in the forbidden locus. Furthermore, for these embeddings, we give a component of the real rank boundary, the hypersurface dividing the minimal typical rank from higher ones. These results generalize to a class of embeddings of $\mathbb P^n\times \mathbb P^1$. Finally, in connection with real rank boundaries, we give a new interpretation of the $2\times n \times n$ hyperdeterminant.
\end{abstract}

\noindent \begin{small}; \end{small}

\section{Introduction}

The study of loci of points in projective space having exceptional properties with respect to a projective variety is a classical theme in projective geometry. Recently the thriving of multilinear algebra in both pure and applied mathematics has revitalized classical geometric questions in the framework of tensor rank. This notion can be stated in full generality for any projective variety. Let $\mathbb K$ be a field of characteristic zero and let $X\subset \mathbb P^N_{\mathbb K}$ be a projective variety. We assume $X$ is non-degenerate, that is, $X$ is not contained in a hyperplane. \\
\indent The $X$-rank of a point $f\in \mathbb P^N_{\mathbb K}$, denoted $\textnormal{rk}_X(f)$, is the minimum integer $s$ such that $f$ is in the span of $s$ distinct points of $X$: 
$$
f \in\langle \ell_1,\ldots, \ell_s \rangle, \mbox{ where }  \ell_i\in X.
$$
\noindent When $X$ is a Veronese or a Segre variety, one retrieves the classical Waring rank of homogeneous polynomials, or symmetric rank of symmetric tensors, and the tensor rank of general tensors. As a consequence, the knowledge of $X$-ranks is of fundamental importance. Indeed, the study of $X$-ranks is motivated by applications to algebraic statistics, complexity theory, quantum information theory, signal processing, and many other fields; see \cite{Lands} and the numerous references therein. \\
\indent Among the possible $X$-ranks appearing in the ambient real projective space of $X$, there are some special ones called {\it typical ranks}, which we define next. 

\begin{definition}[{\bf Typical rank}]

Let $X\subset \mathbb P^N_{\mathbb R}$ be a real projective variety. Let $\mathcal R_s = \lbrace{ f\in \mathbb P^N_{\mathbb R}\ |\ \textnormal{rk}_{X}(f) = s \rbrace}$. These sets are semi-algebraic. If $\mathcal R_s$ contains an open euclidean ball, then $s$ is a {\it typical rank}. 

\end{definition}

In probabilistic terms, a rank is typical when the probability of having a randomly sampled tensor with that rank is positive. Typical ranks naturally arise in applications \cite{tenB}. \\
\indent The analogous notion over $\mathbb K = \mathbb C$, or more generally for an algebraically closed field, yields a unique rank, which is the {\it generic rank}. Over $\mathbb K = \mathbb R$, there might be several typical ranks. An instance of this phenomenon is shown by binary forms \cite{Ble}. There have been several studies concerning typical ranks for Veronese or arbitrary non-degenerate projective varieties; see \cite{BBO, MMSV16}. Even when $X$ is a Veronese variety, this problem is far from being fully solved, and only partial results are known. A challenging question is detecting all typical ranks for a given real projective variety $X$. Whenever a real projective variety has more than one typical rank, the transition between two typical ranks happens across a complex algebraic hypersurface. This hypersurface dividing the smallest typical rank from higher ones is the {\it real rank boundary}. Here is a formal definition. (When $\mathbb K=\mathbb C$ we drop the subscript to denote the projective space.)
 
\begin{definition}[{\bf Real rank boundary}]

Let $X\subset \mathbb P^N_{\mathbb R}$ be a real projective variety. Let $g$ be the complex generic rank of its complexification in $\mathbb P^N$. 
We define $\mathcal{R}_X$ to be the set
$$
\mathcal{R}_X = \lbrace f\in \mathbb P^N_{\mathbb R}\ |\ \textnormal{rk}_{X}(f) = g \rbrace.
$$
\noindent This is a full-dimensional semi-algebraic set. The topological boundary $\partial \mathcal{R}_X$ is the closure of $\mathcal{R}_X$ minus its interior. 
If $X$ has more than one typical rank, $\partial \mathcal{R}_X$ is non-empty and its Zariski closure $\partial_{alg} \mathcal{R}_X$ over $\mathbb C$ is a hypersurface in $\mathbb P^N$. This is the {\it real rank boundary} of $X$. 

\end{definition}

One of the main issues is to identify the irreducible components of this hypersurface. The question of describing this algebraic boundary was addressed and completely solved for binary forms by Lee and Sturmfels \cite[Theorem 4.1]{LS}. \\

\noindent {\bf Numerical ranks and real rank boundaries.}\\
A finer knowledge of the real rank boundary, and more generally, of real rank boundaries dividing higher typical ranks, would provide a more detailed picture of the geometry of ranks of a projective variety. In particular, the numerical estimation of the rank of a tensor would benefit from a better understanding of the semi-algebraic geometry of the real rank boundary. Suppose we have a description of the real rank boundary between two real ranks. Informally, if the equation of this real rank boundary almost vanishes on our tensor, it means that we are in the vicinity of it. Close to the boundary usual numerical methods to compute the rank are {\it less} accurate. \\

\noindent {\bf Aim and structure of the article.}\\
\noindent The aim of the present article is to study some specific loci of bigraded forms and real rank boundaries for surfaces arising from special embeddings of $\mathbb P^1\times \mathbb P^1$. These objects are related by the toric antipolar, see Definition \ref{toric antipolar}. \\
\indent Section 2 recalls the concept of forbidden loci introduced in \cite{CCO}. The study of forbidden loci is motivated by the symmetric Strassen's conjecture \cite[Section 4]{CCO}. In this section, we introduce the {\it Ranestad-Schreyer} loci which relate forbidden loci and {\it cactus decompositions}. We analyze this locus in a specific instance for binary forms in Proposition \ref{binform}, Proposition \ref{tangentToVer}, and Proposition \ref{binform2}. In Section 3, we introduce the toric antipolar construction from the toric catalecticants of \cite{GRV}. We study them in the case of embeddings of $\mathbb P^1\times \mathbb P^1$ given by the $(2,2d)$ line bundles. In the tensor literature, these are known as examples of {\it partially symmetric tensors}; see \cite[Chapter 3]{Lands}. For such embeddings of $\mathbb P^1\times \mathbb P^1$, Theorem \ref{nonred locus and antipolar} shows that the vanishing of the toric antipolar coincides with the Ranestad-Schreyer locus. These are contained in the forbidden locus as proven in Lemma \ref{forbid} and Corollary \ref{RanSch=Forb=Antipol}. These results generalize to the embeddings given by the $(2,2d)$ line bundles on $\mathbb P^n\times \mathbb P^1$ whenever $d$ or $n$ is odd, as observed in Remark \ref{p1xpn}. In Section 4, we discuss typical ranks. Theorem \ref{algbound} gives one component of the real rank boundary. In Section 5, we recall known results on hyperdeterminants and we restate them in our terminology. In Proposition \ref{hypervstang}, we show that the hyperdeterminant of tensors in $\mathbb C^2\otimes \mathbb C^{n}\otimes \mathbb C^n$ is the join of the tangential and the $(n-2)$nd secant variety of the Segre variety. This gives, as far as we know, a new interpretation of these hyperdeterminants. In particular, this gives an alternative description of the real rank boundary between 
the unique two typical ranks in $\mathbb P(\mathbb R^2 \otimes \mathbb R^{n}\otimes \mathbb R^n)$, as remarked in Corollary \ref{newinterpofbound}.

\section{Forbidden loci}

In this section, we introduce some loci of points in projective space and study relations among them. Our first object is the {\it forbidden locus}.

\begin{definition}[{\bf Forbidden locus}]

Let $X$ be a projective variety in $\mathbb P^N$. The {\it forbidden locus} $\mathfrak{F}(f)$ of $f\in \mathbb P^N$ is the subset of points in $X$ 
that are not in the support of any smooth zero-dimensional subscheme $Z\subset X$ of length $\textnormal{rk}_{X}(f)$ with $f\in \langle Z \rangle$. 

\end{definition}
\noindent Equivalently, $X\setminus \mathfrak{F}(f)$ is the set of $\ell$ for which there exists $\lambda \in \mathbb C^{*}$ such that $f+\lambda \ell$ has $X$-rank smaller than $\textnormal{rk}_{X}(f)$. \\
\indent Motivated by Strassen's conjecture \cite[Chapter 5.2]{Lands}, Carlini, Catalisano, and Oneto \cite{CCO} have started the study of forbidden loci in the context of Waring decompositions. This is the situation above for $X$ being the $d$th Veronese of a projective space $\mathbb P^n$ and $N = \binom{n+d}{d}-1$. The present article has been primarily inspired by the example of ternary cubics. Independently in \cite{CCO} and \cite{MMSV16}, the forbidden loci of general ternary cubic forms were determined. \\
\indent One of our objectives is to show that forbidden loci are related to loci that parameterize decompositions with non-reduced structure at a point. These decompositions are examples of cactus decompositions. Such loci were already introduced by Micha\l{}ek and Moon \cite{MM17}, building upon results by Ranestad and Schreyer \cite{RSpolar}. We name them {\it Ranestad-Schreyer loci}. 

\begin{definition}[{\bf Ranestad-Schreyer loci}]\label{defRSlocus}
Let $X$ be a projective variety in  $\mathbb P^N$. The {\it Ranestad-Schreyer locus} of length $s$, denoted $\mathfrak{RS}(f)_s$, is the subset of points $\ell\in X$ with the following property: for every $\ell\in \mathfrak{RS}(f)_s$, there exists a zero-dimensional subscheme $Z\subset X$ of length $s$ with $f\in \langle Z \rangle$, that has $\ell$ as non-reduced point. When we deal with $\mathfrak{RS}(f)_s$, for $s = \textnormal{rk}_{X}(f)$, we drop the subscript and we simply refer to it as Ranestad-Schreyer locus of $f$ omitting the length. 
\end{definition}
\begin{remark}\label{constructible}
The loci above are constructible sets. 
\end{remark}
\indent Equipped with this definition, the forbidden locus of general cubics is described by the following result. 

\begin{proposition}[{\bf \cite{CCO, MMSV16}}]
Let $f$ be a general ternary cubic form. The forbidden locus $\mathfrak{F}(f)$ is closed and contains as irreducible component the Ranestad-Schreyer locus $\mathfrak{RS}(f)$, which is the Cayleyan of $f$. The other component is the dual of the Hessian of $f$.  
\end{proposition}

To interpret Definition \ref{defRSlocus} more geometrically, we recall the notion of open {\it areole} \cite[Section 4.1]{LM17}. The $r$th {\it open areole} $\mathfrak{a}^{\circ}_r(X)$ of $X$ is the union of the spans of smoothable subschemes of $X$ of length at most $r$ supported at some $\ell\in X$. If $f\in {\mathfrak{a}}^{\circ}_r(X) + \sigma^{\circ}_{s-r}(X)$, the join of the open areole and the open secant variety, for some $r\geq 2$, then $\mathfrak{RS}(f)_s$ is non-empty. \\ 
\indent We now consider the situation when $X$ is a rational normal curve of degree $d$. This corresponds to considering binary forms of degree $d$. For the sake of completeness, let us recall the apolar ideal of a form; we refer to \cite[Chapter 1]{IK} for more details. 

\begin{definition}[{\bf Apolar ideal}]
Let $V$ be a complex or real vector space and $f\in S^{d}(V^{*})$, that is, a form of degree $d$. The apolar ideal $f^{\perp}$ is the ideal consisting of all forms in the symmetric algebra $S(V) = \bigoplus_{r\in \mathbb N} S^r(V)$ that annihilate $f$ by differentiation.  
\end{definition}

\begin{proposition}\label{binform}
Let $f$ be a binary form of degree $d$ with maximal complex rank $\textnormal{rk}_{\mathbb C}(f) = d$. Then the forbidden locus $\mathfrak{F}(f)$ and the Ranestad-Schreyer locus of length two $\mathfrak{RS}(f)_2$ coincide. Whenever $f$ does not have the maximal complex rank, the Ranestad-Schreyer locus $\mathfrak{RS}(f)_2$ is empty and the forbidden locus is not. 
\begin{proof}
From Definition \ref{defRSlocus}, it follows that $\mathfrak{RS}(f)_2\neq \emptyset$ if and only if $f$ is spanned by a double point if and only if $f\in \tau(\nu_d(\mathbb P^1))\setminus \nu_d(\mathbb P^1)$, where $\tau(\nu_d(\mathbb P^1))$ is the tangential variety of the rational normal curve. Binary forms in $\tau(\nu_d(\mathbb P^1))\setminus \nu_d(\mathbb P^1)$ are the only ones that do have complex rank $d$; see, for example, \cite[Proposition 19]{BHMT}. If $f$ has maximal complex rank, \cite[Theorem 3.5]{CCO} shows that $\mathfrak{F}(f)$ is the {\it unique} point of tangency of $f$. (The case $d=2$ is special: forbidden locus and Ranestad-Schreyer locus of $f$ are both given by the {\it two} points of tangency of the two tangent lines to the conic $\nu_2(\mathbb P^1)\subset \mathbb P^2$ passing through $f$.) This gives the desired equality. Whenever $f$ is not of maximal complex rank, $f$ is not in $\tau(\nu_d(\mathbb P^1))\setminus \nu_d(\mathbb P^1)$ and hence $\mathfrak{RS}(f)_2 = \emptyset$. In this case, the forbidden locus is not empty by the case analysis in \cite[Theorem 3.5]{CCO}.
\end{proof}
\end{proposition}

A consequence of Proposition \ref{binform} is the following. 

\begin{proposition}\label{tangentToVer}

Let $X$ be a Veronese variety of degree $d$ and suppose $f\in \tau(X)\setminus X$. Then $\mathfrak{RS}(f)_2 \subseteq \mathfrak{F}(f)$. 
\begin{proof}
From Definition \ref{defRSlocus}, $\ell\in \mathfrak{RS}(f)_2$ if and only if $f$ is in the tangent space of $X$ at $\ell$. For any such $\ell$, consider 
a suitable rational normal curve $\nu_d(\mathbb P^1)\subset X$ passing through $\ell$ and such that $f$ is in the tangential of it. Hence $f$ can be regarded as a binary form $\tilde f$ of degree $d$ of maximal complex rank. By definition of $\tilde f$ and Proposition \ref{binform}, $\ell \in \mathfrak{RS}(\tilde f)_2=\mathfrak{F}(\tilde f)$, where the latter is contained in $\mathfrak{F}(f)$. Thus $\mathfrak{RS}(f)_2 \subseteq \mathfrak{F}(f)$. 
\end{proof} 

\end{proposition}

\begin{proposition}\label{binform2}
Let $d\geq 3$. Let $f$ be a binary form of degree $d$ with complex rank $\textnormal{rk}_{\mathbb C}(f) = d$. Then $\mathfrak{RS}(f)=\nu_d(\mathbb P^1)$. 

\begin{proof}
Up to change of coordinates we may assume that $f = x^{d-1}y$; this is tangent to $[x^d]$. The apolar ideal of $f$ is $f^{\perp} = (x^d,y^2)$. A point $\ell$ is in $\mathfrak{RS}(f)$ if and only if there exists a homogeneous polynomial of degree $d$ in $f^{\perp}$ such that it has the linear form dual to $\ell$ as a factor with multiplicity strictly larger than one. This imposes linear conditions on the coefficients of a generic polynomial in $f^{\perp}_d$, which have a solution for every point in $\nu_d(\mathbb P^1)$. 
\end{proof}
\end{proposition}

We shift gears to concrete examples of forbidden loci. As observed in Remark \ref{constructible}, forbidden loci are constructible sets. In the special situation where the closure of the locus of forms with higher rank than the generic complex rank is given by the closure of a single $\textnormal{GL}$-orbit, the forbidden locus is contained in a proper Zariski closed subset:  

\begin{proposition}\label{glorbits}

Let $X$ be a Veronese variety. Let $m$ be an integer strictly larger than the generic complex rank $g$ and let $W_m$ be the closure of the locus of forms of rank $m$. Suppose that $W_m$ is the closure of the $\textnormal{GL}$-orbit of a form $f$.  Then the forbidden locus of $f$ is contained in a proper closed subset of $X$. 

\begin{proof}
By \cite[Theorem 7]{BHMT}, the join $W_m+X$ is contained in $W_{m-1}$ for all $m>g$. Hence the cone $f+X$ sits inside $W_{m-1}$. The closure of the $\textnormal{GL}$-orbit of this cone coincides with $W_m+X$. Thus, the general point $f+\lambda \ell$ has complex rank $m-1$. This implies that the forbidden locus $\mathfrak{F}(f)$ is contained in a proper closed subset of $X$. 
\end{proof}

\end{proposition}

Two examples where the forbidden locus is closed are ternary and quaternary cubics of maximal rank.

\begin{example}[{\bf Cubics of maximal rank}]
Ternary cubics whose components are a smooth conic and a tangent line are the only cubics which have maximal complex rank five. Hence $W_5$ is the $\textnormal{GL}$-orbit of such a cubic. In  \cite[Theorem 3.18]{CCO}, it is shown that the forbidden locus for $f$ in this orbit is one point. Similarly, cubics in $\mathbb P^3$ whose components are a smooth quadric and a tangent plane are the only cubics which have maximal complex rank seven. Hence $W_7$ is the $\textnormal{GL}$-orbit of such a cubic. The forbidden locus of these cubics is a single point. Indeed, by the same argument in the proof of \cite[Theorem 3.18]{CCO}, based on the second Bertini's theorem, the forbidden locus $\mathfrak{F}(f)$ is the point dual to this tangent plane. 
\end{example}

We now explicitly compute forbidden loci in the next two examples of quartics of maximal complex rank seven. We rely on the classification of complex ranks of ternary quartics by Kleppe \cite[Chapter 3]{Kleppe}. Before we proceed, we briefly recall the construction of the classical catalecticant; we refer to \cite[Chapter 3.5]{Lands} for more details. 

\begin{definition}[{\bf Classical catalecticant}]\label{classicalcat}
Let $V$ be a complex or real vector space and $f\in S^{d}(V^{*})$, that is, a form of degree $d$. For every $1\leq r \leq d$, we consider the linear map 
$$
C_{f, r}\in S^r(V^{*})\otimes S^{d-r}(V^{*})
$$
\noindent induced by $f$ and defined by $C_{f, r}(g) = g(f)$. This is called catalecticant. When $d = 2k$, the catalecticant $C_{f,k}$ is called middle catalecticant and denoted by $C_f$. 
\end{definition}

\begin{example}[{\bf A reducible quartic}]\label{redquartic}
The forbidden locus of the reducible ternary quartic $f = y^2(xy+z^2)$ of maximal complex rank is one point. Indeed, its middle catalecticant $C_f$ has rank three. To determine the forbidden locus, let us consider an arbitrary linear form $\ell = (ax+by+cz)$, which is identified with the projective point $[a:b:c]$. We compute the middle catalecticant $C_{f + \lambda \ell^4}$ of the form $f+\lambda \ell^4$, with $\lambda\in \mathbb C^{*}$. This matrix has rank at most four for any choice of $a,b,c,\lambda$. Moreover, \cite[Theorem 3.2]{Kleppe} states that whenever the middle catalecticant of a quartic $g$ has a two-dimensional kernel, $g$ has either rank four or six. By this result, if the rank of $C_{f + \lambda \ell^4}$ is four, then $f+\lambda \ell^4$ has complex rank six. (Indeed, it cannot be less, because it would contradict the fact that the complex rank of $f$ is seven.) Its rank is three only when $a=c=0$. We check with \textnormal{\texttt{Macaulay2}} \cite{M2} that there exist projective transformations sending $f$ to $f+\lambda y^4$ for every $\lambda\in \mathbb C^{*}$. Hence the complex rank of $f+\lambda y^4$ is seven. Thus the forbidden locus is $[0:1:0]$. \\
\end{example}

\begin{example}[{\bf Yet another quartic}]
The forbidden locus of the ternary quartic $f = (x+y)^4+(x^3+y^3)z$ of maximal complex rank is two points. Indeed, its middle catalecticant $C_f$ has rank five. As above, we consider $f+\lambda\ell^4$ and its middle catalecticant $C_{f + \lambda \ell^4}$. We have $\det (C_{f + \lambda \ell^4}) = -746496\lambda c^4$. If $c\neq 0$, then $C_{f + \lambda \ell^4}$ has full rank. Moreover, \cite[Theorem 3.7]{Kleppe} states that whenever the middle catalecticant of a quartic $g$ has zero-dimensional kernel, $g$ has rank six. This result implies that $f + \lambda \ell^4$ has rank six. This means that $\ell$ appears in a minimal decomposition of $f$. If $c=0$ and $a,b\neq 0$, then there exists a unique nonzero value $\tilde \lambda$ such that $C_{f + \tilde{\lambda} \ell^4}$ has rank four. This implies that $f + \tilde\lambda \ell^4$ has rank six by \cite[Theorem 3.2]{Kleppe}; see Example \ref{redquartic} for an analogous application of this result. If $c=0$ and $a=0$ or $b=0$, then $C_{f + \lambda \ell^4}$ has rank five. Its kernel is generated by $z^2$. Now, \cite[Theorem 3.6]{Kleppe} states that whenever the middle catalecticant of a quartic $g$ has kernel generated by a double line, $g$ has rank seven. This theorem shows that $f + \lambda \ell^4$ has complex rank seven. This means that the only points in the forbidden locus are $[1:0:0]$ and $[0:1:0]$. 
\end{example}

\section{Antipolars and toric apolarity}

In this section, we introduce antipolars in the context of toric apolarity. The aim is to deduce a further connection between forbidden loci and Ranestad-Schreyer loci, and new results on typical ranks for toric varieties in Section \ref{typtoric}. \\ 
\indent We now give the classical definition of antipolar of a form $f$ of even degree. 

\begin{definition}[{\bf Classical antipolar \cite{BBO, MMSV16}}]

Let $f\in S^{2k}(V^{*})$ and $C_f$ be its middle catalecticant. Assume that $C_f$ is an isomorphism. 
The antipolar of $f$ is defined as 
$$
\Omega(f)(\ell) = \det(C_{f+\ell^{2k}}) - \det(C_f).
$$ 
\noindent The forms $\Omega(f)$ and $f$ have the same degree. 

\end{definition}

\begin{remark}\label{classident}
We have the classical identity 
$$
v^{t} A^{-1} u = \frac{1}{\det(A)}\left[\det( A + uv^{t}) - \det(A)\right],
$$
\noindent where $A$ is an invertible $k\times k$ complex (or real) matrix and $u,v$ are vectors in $\mathbb C^k$ (or in $\mathbb R^k$). Then the form 
$\Omega(f)$ above is defined up to scaling as the bilinear map given by the inverse of $C_f$. 
\end{remark}
\noindent This construction generalizes the one of the dual quadric. \\
\indent In \cite{MM17}, Micha\l{}ek and Moon generalized a result of Ranestad and Schreyer \cite[Lemma 2.3]{RSpolar}, showing that the antipolar of a form $f\in S^{2d}(V^{*})$, gives the full knowledge of a non-reduced structure at a point of a zero-dimensional scheme spanning $f$. In our terminology, they showed that, under special assumptions, the Ranestad-Schreyer locus is given by the vanishing of the antipolar. 

\begin{proposition}[{\bf\cite[Proposition 2.6]{MM17}}]\label{Nonreducedlocus}

Let $d\in \lbrace 1,2,3,4 \rbrace$ and $f\in S^{2d}(V^{*})$  be a ternary form of degree $2d$. An apolar scheme $Z$ to $f$ of degree $\crk(f)$ has non-reduced structure at a point $\ell \in Z$ if and only if $\Omega(f)(\ell) = 0$. Equivalently, the Ranestad-Schreyer locus $\mathfrak{RS}(f)$ coincides with $\lbrace \Omega(f)(\ell) = 0\rbrace$. 

\end{proposition}

\indent Our objective is to give a more general version of Proposition \ref{Nonreducedlocus} for the apolarity of toric varieties.  We follow the approach by Gallet, Ranestad, and Villamizar \cite{GRV}, which we fully recall here. \\
\indent Let $X$ be a toric variety and let $C(X)$ be the Cox ring of $X$. This is the graded $\mathbb C$-algebra of sections of all line bundles on $X$. The ring $C(X)$ is naturally graded by the Picard group $\textnormal{Pic}(X)$ of $X$, the group of isomorphism classes of line bundles on $X$. The Cox ring of a toric variety is a polynomial ring. \\
\indent Let $S=C(X)$ and let $T$ be its {\it dual ring}, that is, for each $A\in \textnormal{Pic}(X)$, $T_A$ is the dual vector space of $S_A$. For $f\in S_A$, the corresponding hyperplane in $T_A$ is denoted by $H_f$. 

\begin{definition}[{\bf \cite[Definition 1.2]{GRV}}]\label{apolardef}

A subscheme $Z\subset X$ is {\it apolar} to $f\in S_A$ if and only if the ideal $I_{Z,A} = \lbrace g\in T_A | g(Z) = 0\rbrace$ of $T_A$ is contained in $H_f$. 

\end{definition}

To fully recover the classical setting, the next step is having an analogue of the apolarity lemma \cite[Lemma 1.15]{IK}. To this aim, let us introduce a partial ordering on $\textnormal{Pic}(X)$. For $A,B\in \textnormal{Pic}(X)$, we let $A > B$ if the line bundle $A-B$ has global sections. Equipped with this ordering, for $f\in S_A$, let us define the ideal $I_f$ by describing its graded components: 
$$
I_{f,B} = 
  \begin{cases}
  H_f:T_{A-B}=\lbrace g\in T_B \ | \ gT_{A-B}\subseteq H_f\rbrace & \quad \text{if } A > B,\\
  T_B & \quad \text{otherwise}.\\
  \end{cases}
$$

\noindent Thus we define $I_f = \oplus_{B\in \textnormal{Pic}(X)} I_{f,B}$, which is an ideal of $T$. Analogously, for a given subscheme $Z\subset X$, $I_Z =  \oplus_{B\in \textnormal{Pic}(X)} I_{Z,B}$. The toric apolarity lemma is the following. 

\begin{lemma}[{\bf\cite[Lemma 1.3]{GRV}}]\label{apolarlemma}
If $Z\subset X$ and $f\in S_A$, then $I_Z \subset I_f$ if and only if $I_{Z,A}\subseteq I_{f,A} = H_f$. 
\end{lemma}
Given a form $f\in S_A$, we introduce a family of linear maps that replace catalecticant matrices in this setting. For any pair of classes 
$A,B\in \textnormal{Pic}(X)$, one defines a linear map 
$$
\phi_{f,B}: T_B \rightarrow S_{A-B},\\
g\mapsto g(f),
$$
\noindent such that $g(f)(h) = hg(f)\in \mathbb C$, for $h\in T_{A-B}$. 

\begin{remark}
The notation of these catalecticant matrices differs from the one of Definition \ref{classicalcat} in order to maintain the usual notation. To recover Definition \ref{classicalcat}, put $S_{A-B} = S^{d-r}(V^{*})$ and $T_B = S^{r}(V)$. 
\end{remark}

\begin{remark}
From the definitions above, it readily follows that $\ker \phi_{f,B} = I_{f,B}$ for every $B\in \textnormal{Pic}(X)$.  
\end{remark}

\begin{remark}[{\bf Non-abelian apolarity}]
Toric apolarity might be seen as particular case of a more general vector bundle construction \cite{LO11}. Let $X$ be a projective variety and $L$ a very ample line bundle, which gives an embedding $X\subset \mathbb P(H^0(X,L)^{*})=\mathbb P(W)$. Let $E$ be a vector bundle on $X\subset \mathbb P(W)$. For a given $f\in W$, we can construct a linear map $A_f$. This construction is linear in $f\in W$ and it is as follows. We first consider the natural contraction map on global sections induced by $E\otimes E^{*}\otimes L \rightarrow L$, $H^0(E)\otimes H^0(E^{*}\otimes L) \rightarrow H^0(L)$. Viewing this linear map as a tensor and taking another flattening of it, we obtain a linear map: 
$$
H^0(E) \otimes H^0(L)^{*} \rightarrow H^0(E^{*}\otimes L)^{*}.
$$

\noindent This gives a linear map

$$
A_f: H^{0}(E) \rightarrow H^{0}(E^{*}\otimes L)^{*},
$$

\noindent which depends linearly on $f\in H^{0}(L)^{*}$. The map $A_f$ is the catalecticant. In this framework, called {\it non-abelian apolarity}, an analogous apolarity lemma holds; see \cite[Proposition 5.4.1]{LO11}. The classical apolarity is recovered for $X$ being the projective space and $L, E$ being line bundles. Our use of toric apolarity is obtained from this general construction: for $X$ being a projective toric variety along with $L, E$ being very ample line bundles on $X$. 
\end{remark}

\begin{definition}
Let $B\in \textnormal{Pic}(X)$ be a very ample line bundle. We denote by $\mathfrak{b}: X\hookrightarrow \mathbb P(H^0(X,B)^{*})$ the corresponding embedding. Note that the embedding corresponding to $2B\in \textnormal{Pic}(X)$ is the composition of $\mathfrak{b}$ and the second Veronese $\nu_2$. \end{definition}

We introduce antipolars in the context of toric apolarity. 

\begin{definition}[{\bf Toric antipolar}]\label{toric antipolar}

We keep the notation from above. Let $f\in S_{2B}$ be a form such that $\phi_{f,B}$ is an isomorphism. Then the antipolar of $f$ is
$$
\Omega(f)(\ell) = \det(\phi_{f + \nu_2(\mathfrak{b}(\ell)),B}) -\det(\phi_{f,B}),
$$
\noindent where $\ell \in X$. 

\end{definition}

\begin{remark}

Let $f\in S_{2B}$ be a form such that $\phi_{f,B}$ is an isomorphism. The antipolar $\Omega(f)$ is a form in $S_{2B}$ by Remark \ref{classident}.

\end{remark}

For the sake of simplicity of notation, for a given $f\in S_{2B}$, we refer to its complex $\nu_2(\mathfrak{b}(X))$-rank as $X$-rank of $f$. A general form $f\in S_{2B}$ has a full rank catalecticant $\phi_{f,B}$.  

\begin{theorem}\label{nonred locus and antipolar}

Let $B\in \textnormal{Pic}(X)$ be a very ample line bundle.  Let $f\in S_{2B}$ be a general form and suppose that  its complex $X$-rank is $\dim_{\mathbb C} T_B$, the size of $\phi_{f,B}$. Let $Z\subset X$ be a zero-dimensional subscheme of $X$, whose degree is equal to the dimension of $T_B$, and apolar to $f$. Then $Z$ has a non-reduced structure at $\ell \in Z$ if and only if $\Omega(f)(\ell) = 0$. In other words, the Ranestad-Schreyer locus $\mathfrak{RS}(f) = \lbrace \Omega(f)(\ell) = 0 \rbrace$.

\begin{proof}

Let $\tilde{Z}$ be the scheme defined by $I_{\tilde{Z}} = I_{Z}: I_{\ell}$, where $\ell\in Z$. Since the linear map $\phi_{f,B}$ is an isomorphism by assumption, we have
$\ker \phi_{f,B} =0$. The degree of $Z$ coincides with the dimension of $T_B$; thus the degree of $\tilde{Z}$ is strictly smaller than this dimension and hence we have 
a section $g\in T_B$ that vanishes at $\tilde{Z}$. This implies that $g I_{\ell} \subset I_Z$ and then $gI_{\ell, B}\subset I_{Z}$. Since $Z$ is apolar to $f$, by Definition \ref{apolardef}, $I_{Z,2B}\subseteq I_{f,2B}$. Thus $g I_{\ell, B} \subset I_{f,2B} = H_f$. This means that $g (\phi_{f,B}(I_{\ell, B})) = 0$. This defines the zeros of the section $g$, since $I_{\ell,B}$ is a hyperplane in $T_B$. Consequently, $Z$ is reduced at $\ell$ if and only if $\mathfrak{b}(\ell)\notin \phi_{f,B}(I_{\ell, B})$ if and only if $\phi_{f,B}^{-1}(\mathfrak{b}(\ell)) \notin I_{\ell,B}$. Thus $Z$ has a non-reduced structure at $\ell$ if and only if the pairing $\langle \mathfrak{b}(\ell), \phi_{f,B}^{-1}(\mathfrak{b}(\ell)) \rangle=\Omega(f)(\ell) = 0$. 
\end{proof}
\end{theorem}

The next lemmas describe part of the forbidden locus in the situation of Theorem \ref{nonred locus and antipolar}.

\begin{lemma}\label{WaringLocus}

Let $f\in S_{2B}$ be a general form and suppose that  its complex $X$-rank is $\dim_{\mathbb C} T_B$, the size of $\phi_{f,B}$. Let $\ell\in X$. If $\ell \notin \mathfrak{F}(f)$ then there exists $\lambda \in \mathbb C^*$ such that $\det(\phi_{f+\lambda \nu_2(\mathfrak{b}(\ell)),B}) = 0$. 

\begin{proof}

Suppose that $\ell \notin \mathfrak{F}(f)$. Then there exists a minimal smooth apolar scheme to $f$, where $\nu_2(\mathfrak{b}(\ell))$ appears. This means $f =\sum_{i=1}^{\dim_{\mathbb C} T_B} \lambda_i \nu_2(\mathfrak{b}(\ell_i))$, where there exists $1\leq j\leq \dim_{\mathbb C} T_B$ such that $\ell = \ell_j$. Choosing $\lambda = -\lambda_j\in \mathbb C^*$, we have that $f+\lambda \nu_2(\mathfrak{b}(\ell))$ has smaller rank than the size of $\phi_{f,B}$. Hence $\phi_{f+\lambda \nu_2(\mathfrak{b}(\ell)),B}$ is degenerate.  
\end{proof}
\end{lemma}

\begin{lemma}\label{forbid}

Let $f\in S_{2B}$ be a general form and suppose that  its complex $X$-rank is $\dim_{\mathbb C} T_B$, the size of $\phi_{f,B}$. Then the forbidden locus $\mathfrak{F}(f)$ contains the locus of the antipolar $\Omega(f)$. 

\begin{proof}
Let $\ell\notin \mathfrak{F}(f)$. By Lemma \ref{WaringLocus}, there exists $\lambda \in \mathbb C^*$ such that $\det(\phi_{f + \lambda \nu_2(\mathfrak{b}(\ell)),B})=0$. Note that $\lambda\Omega(f)(\ell) = \det(\phi_{f + \lambda \nu_2(\mathfrak{b}(\ell)),B}) - \det(\phi_{f,B})$. Hence $\lambda \Omega(f)(\ell) = -\det(\phi_{f,B})\neq 0$, which implies $\Omega(f)(\ell)\neq 0$. Thus $\lbrace \Omega(f)(\ell) = 0\rbrace \subseteq \mathfrak{F}(f)$.  
\end{proof}
\end{lemma}

As direct corollary, we have the following fact. 

\begin{corollary}\label{RanSch=Forb=Antipol}
Let $B\in \textnormal{Pic}(X)$ be a very ample line bundle.  Let $f\in S_{2B}$ be a general form and suppose that its complex $X$-rank is $\dim_{\mathbb C} T_B$, the size of $\phi_{f,B}$. The Ranestad-Schreyer locus $\mathfrak{RS}(f) = \lbrace \Omega(f)(\ell) = 0\rbrace$ is contained in the forbidden locus $\mathfrak{F}(f)$. In particular, the forbidden locus is non-empty. 
\end{corollary}

In other words, in the situation of Corollary \ref{RanSch=Forb=Antipol}, points that belong to cactus decompositions of $f$, in which they appear with a non-reduced structure, are also forbidden to appear in any minimal smooth decomposition of $f$.  

\begin{remark}
By a result of Clebsch \cite[Corollary 3.5.1.5]{Lands}, the generic rank of forms $f\in S^{4}(V^{*})$, with $3\leq \dim(V^{*})\leq 5$, is $\binom{\dim(V^{*})+1}{2}$, where the latter is the dimension of $S^{2}(V^{*})$. In these cases, Corollary \ref{RanSch=Forb=Antipol} holds as well. 
\end{remark}

\section{Typical ranks of embeddings of $\mathbb P^1 \times \mathbb P^1$}\label{typtoric}

We study a specific instance of toric apolarity and we give an application to typical ranks. Throughout this section, we let $X = \mathbb P^1 \times \mathbb P^1$, the Segre surface. The Picard group of $X$ is isomorphic to $\mathbb Z^2$. Its Cox ring $C(X)$ is a polynomial ring in four variables $x,y,z,w$, equipped with the bigrading $\deg(x) = \deg(y) = (1,0)$ and $\deg(z) = \deg(w) = (0,1)$. As above, the ring $C(X)$ can be written as a direct sum of its bihomogeneous components, that is $C(X) = \bigoplus_{A = (u,v)\in \mathbb Z^2} C(X)_A$, where $C(X)_A$ is of bidegree $A=(u,v)$. We set $S = \mathbb C[x,y,z,w]$ and $T = \mathbb C[\partial_x, \partial_y, \partial_z, \partial_w]$. 
\noindent Over $X$, every divisor of class $A = (u,v)$ with $u,v > 0$ is very ample and determines the {\it Segre-Veronese embedding}: 
$$
\nu_{(u,v)}: \mathbb P^1 \times \mathbb P^1 \rightarrow \mathbb P^{uv+u+v}=\mathbb P(C(X)_A),
$$
$$
(l_1,l_2) \mapsto l_1^ul_2^v,
$$
\noindent where $l_1, l_2$ are points in the two copies of $\mathbb P^1$ respectively. The points in $\mathbb P(C(X)_A)$ are examples of {\it partially symmetric tensors} \cite[Chapter 3]{Lands}.

\begin{example}
Let $f = 4x^2z^2+6x^2zw+2x^2w^2+8xyz^2+7xyzw+5xyw^2+3y^2z^2+7y^2zw+2y^2w^2 \in S_{2B}$, where $B=(1,1)$.  The catalecticant $\phi_{f,B}: T_{B} \rightarrow S_B$ is given by the $4\times 4$ symmetric matrix: 
$$
\phi_{f,B} = \begin{pmatrix}
16 & 12 & 16 &  7 \\ 
12 & 8 & 7 & 10 \\
16 & 7 & 12 & 14 \\
7 & 10 & 14 & 8  \\
\end{pmatrix}. 
$$
\noindent We now introduce unknown coefficients $s_1,s_2,t_1, t_2$. Let $l_1 = s_1x+s_2y$ be a linear form of bidegree $(1,0)$ and $l_2 = t_1z+t_2w$ a linear form of bidegree $(0,1)$. The antipolar of $f$ at $\nu_2(\mathfrak{b}(\ell)) = l_1^2l_2^2$ is the following form of bidegree $2B=(2,2)$, in the bigraded variables $s_1,s_2,t_1,t_2$: 
$$
\Omega(f)(\ell) = -1728s_1^2t_1^2-1104s_1s_2t_1^2-1760s_2^2t_1^2 
+80s_1^2t_1t_2+12272s_1s_2t_1t_2-2144s_2^2t_1t_2
$$
$$
-4400s_1^2t_2^2-2048s_1s_2t_2^2-1344s_2^2t_2^2.
$$
\end{example}

The dimensions of secant varieties of Segre-Veronese embeddings of the Segre surface $X$ were determined by 
Catalisano, Geramita, and Gimigliano \cite[Corollary 2.3]{CGG}. 

\begin{theorem}[{\bf Catalisano, Geramita, and Gimigliano}]

Let $f\in S_A$ be a general form. Then the generic complex $X$-rank is

$$
\textnormal{rk}_X(f) = 
  \begin{cases}
  \big \lceil{ \frac{(u+1)(v+1)}{3}}\big\rceil & \quad \text{if } A= (u,v) \neq (2,2d),\\
    2d+2  & \quad \text{if } A = (2,2d).\\
  \end{cases}
$$

\end{theorem}

\begin{remark}
Note that when $A = 2B$, with $B = (1,d)$, the complex $X$-rank of the general form $f\in S_{2B}$ is the same 
as the size of $\phi_{f,B}$. Thus Corollary \ref{RanSch=Forb=Antipol} holds. 
\end{remark}

\begin{remark}\label{mindeg}
In the situation of the previous remark, that is, when $B=(1,d)$, the varieties $\mathfrak{b}(X) \subset \mathbb P^{2d+1}$ are remarkable surfaces. They are 
rational normal scrolls of minimal degree $2d = \textnormal{codim}(\mathfrak{b}(X))+1$. As they are toric varieties, the degree can be easily seen from the number of normalized volume triangles that sit in their corresponding polytopes, which are rectangles of size $1\times d$. 
\end{remark}

Minimal degree varieties are classical varieties that were classified by Bertini and Del Pezzo. More recently, Blekherman, Smith, and Velasco \cite[Theorem 1.1]{BSV} found this beautiful characterization from the real algebraic geometry perspective:

\begin{theorem}[{\bf Blekherman, Smith, and Velasco}]\label{bsv}
Let $X\subset \mathbb P^N_{\mathbb R}$ be a real irreducible non-degenerate projective variety such that the set of real points is Zariski dense. Every non-negative real quadratic form on $X(\mathbb R)$ is a sum of squares of linear forms if and only if $X$ is a variety of minimal degree.
\end{theorem}

\begin{remark}
Let $X = \mathbb P^1_{\mathbb R}\times \mathbb P^1_{\mathbb R}$ and let $B$ be a line bundle. Note that forms of bidegree $2B$ on $X$ are quadratic forms on the varieties $\mathfrak{b}(X)$. For $B=(1,d)$, Remark \ref{mindeg} and Theorem \ref{bsv} (see \cite[Example 5.6]{BSV} for more details) imply that the cone $P_{2B}$ of non-negative forms of bidegree $2B$ on $X(\mathbb R)$ coincides with the cone of sum of squares $\Sigma_{2B}$. In particular, the Segre-Veronese orbitope $P_{2B}^{\vee}$, which is the cone spanned by $\nu_2(\mathfrak{b}(\ell))$ for $\ell \in X(\mathbb R)$, coincides with the dual cone $\Sigma_{2B}^{\vee}$, consisting of forms $f$ of bidegree $2B$ whose middle toric catalecticant $\phi_{f,B}$ is positive semi-definite. 
\end{remark}

We work over $\mathbb R$ and discuss typical ranks with respect to embeddings of $X = \mathbb P^1_{\mathbb R}\times \mathbb P^1_{\mathbb R}$. Henceforth we assume $B=(1,d)$. For a real form $f\in S_{2B}$ we refer to its real rank with respect to the real variety $\nu_2(\mathfrak{b}(X))$ as its real $X$-rank. We first recall a result by Reznick \cite[Theorem 4.6]{Rez}, which describes the relationship between the real rank of a form $f$ and its middle catalecticant $C_f$. This result generalizes to our context since $P_{2B}^{\vee} = \Sigma_{2B}^{\vee}$. We state it in our terminology. 

\begin{theorem}[{\bf Reznick}]\label{rezth}
Let $f\in S_{2B}$ be a real form. Suppose that its middle toric catalecticant $\phi_{f,B}$ is positive semi-definite. Then the real $X$-rank of $f$ coincides with the rank of $\phi_{f,B}$. 
\end{theorem}

We establish the analogue of \cite[Lemma 6.4]{MMSV16} for toric apolarity.

\begin{proposition}\label{changeofsignAntip}
Let $f\in S_{2B}$ be a general real form and suppose that its complex $X$-rank is $\dim_{\mathbb C} T_B$. Moreover, assume that $f$ is not in $P_{2B}^{\vee}$. If $\lbrace \Omega(f)(\ell)=0\rbrace = \emptyset$ over $\mathbb R$ then the real $X$-rank of $f$ is strictly larger than $2+2d$. Furthermore, if $\lbrace \Omega(f)(\ell)=0\rbrace \neq \emptyset$ over $\mathbb R$ and the toric catalecticant $\phi_{f,B}$ has signature $(1+2d,1)$, then the real $X$-rank of $f$ is $2+2d$. 
\begin{proof}
Suppose the real $X$-rank of $f$ is $2+2d$. Since $f$ is not in $P_{2B}^{\vee} = \Sigma_{2B}^{\vee}$, there exists a real point $\ell \in X(\mathbb R)$ such that $\phi_{f+\nu_2(\mathfrak{b}(\ell)),B}$ is degenerate. Hence $\Omega(f)(\ell) = - \det \phi_{f,B}$. Moreover, there exists a real point $\ell'$ such that $\phi_{-f+\nu_2(\mathfrak{b}(\ell')),B}$ is degenerate. Note that we have $\Omega(f)(\ell') = -\Omega(-f)(\ell') = \det \phi_{f,B}$. The first equality follows from the fact that $\Omega(f)(\ell)$ is of total degree $2d+1$ in the coefficients of $f$. Thus $\lbrace \Omega(f)(\ell)=0\rbrace \neq \emptyset$.\\
\indent For the second statement, if $\lbrace \Omega(f)(\ell)=0\rbrace \neq \emptyset$, there exists a real $\ell'\in X(\mathbb R)$ such that $\det \phi_{f+\nu_2(\mathfrak{b}(\ell')), B}=0$. Hence the addition of the rank-one matrix $\phi_{\nu_2(\mathfrak{b}(\ell')),B}$ makes the original matrix $\phi_{f,B}$ positive semi-definite. Theorem \ref{rezth} implies that $f+\nu_2(\mathfrak{b}(\ell'))$ has real $X$-rank equal to the matrix rank of $\phi_{f+\nu_2(\mathfrak{b}(\ell')),B}$, which is $1+2d$. Thus the real $X$-rank of $f$ is exactly $2+2d$ and not less, as its complex $X$-rank is $2+2d$. 
\end{proof}

\end{proposition}

\begin{theorem}\label{algbound}

For every $d\geq 1$ and $B=(1,d)$, the real rank boundary of the real variety $\nu_2(\mathfrak{b}(X))$ is non-empty. One 
of its components is the discriminant of the antipolar $\Omega(f)$.  

\begin{proof}
The discriminant of the form $\Omega(f)$ of bidegree $2B = (2,2d)$ is the boundary between non-negative antipolars on $X(\mathbb R)$, and antipolars that have change of sign. By Proposition \ref{changeofsignAntip}, this also divides forms (whose $\phi_{f,B}$ has signature $(1+2d,1)$) of real $X$-rank $2+2d$, from forms (whose $\phi_{f,B}$ has signature $(1+2d,1)$) of higher real $X$-rank. This means that such a discriminant is a component of the real rank boundary of $\nu_2(\mathfrak{b}(X))$. 
\end{proof}
\end{theorem}
 
\begin{remark}
The existence of a form $f$ whose toric catalecticant $\phi_{f,B}$ has signature $(1+2d,1)$ of real $X$-rank $3+2d$ and the fact that $3+2d$ is a typical rank can be proven analogously as in the proof of \cite[Proposition 4.4]{BBO}. The idea is to produce a catalecticant $\phi_{f,B}$ with signature $(1+2d,1)$ which cannot be updated to a positive semi-definite matrix by the addition of an arbitrary rank-one matrix. Here we sketch how this can be achieved. In the notation adopted at the beginning of this section, let us fix the basis of $S_{B}$ to be $\lbrace xz^d, yz^d, xz^{d-1}w, yz^{d-1}w, \ldots, yw^d\rbrace$, that is, it is lexicographic in $z$ and $w$. Let us consider a form $f = e_{d,0}x^2z^{2d}+f_{d,0}y^2z^{2d}+c_{d,0}xyz^{2d}+e_{d,1}x^2z^{2d-2}w^2+f_{d,1}y^2z^{2d-2}w^2+\cdots+e_{0,d}x^2w^{2d}+f_{0,d}y^2w^{2d}$, where the only {\it mixed} term in $x,y$ is $c_{d,0}xyz^{2d}$. Scalar multiples of $e_{d,0},f_{d,0},c_{d,0}$ appear in the upper-left $2\times 2$-minor of $\phi_{f,B}$ and they can be chosen as in \cite[Proposition 4.4]{BBO}, that is, in such a way that the minor has a positive and a negative eigenvalue and that it cannot be non-negative after updating it with a rank-one matrix. The rest of the coefficients is chosen so that $e_{i,j}$ and $f_{i,j}$ are positive and much smaller than $e_{i-1,j+1}$ and $f_{i-1,j+1}$ respectively. This gives a matrix with signature $(1+2d,1)$, since the signature is given as sign changes of the determinants of the principal minors. (Because of the choice of the coefficients and the combinatorial structure of $\phi_{f,B}$, there cannot be sign changes, except the only one in the first upper-left $2\times 2$-minor chosen above.) Now, the same arguments in the proof of \cite[Proposition 4.4]{BBO} show that such a form has real $X$-rank $3+2d$ and that this rank is typical. 
\end{remark}

In terms of typical ranks of the real projective varieties $\nu_2(\mathfrak{b}(X))$ we obtain:

\begin{corollary}

For every $d\geq 1$ and $B=(1,d)$, there exists more than one typical rank with respect to the real varieties $\nu_2(\mathfrak{b}(X))$. 

\end{corollary}

\begin{remark}\label{p1xpn}
Abrescia \cite[Remark 3.9]{Abrescia} showed that the smallest secant variety of the embedding of $\mathbb P^n\times \mathbb P^1$ with the very ample
line bundle $2B = (2,2d)$ filling the ambient space $\mathbb P^{\frac{(2d+1)(n+1)(n+2)}{2}-1}$ is the $(d+1)(n+1)$th secant. The variety $\mathbb P^n\times \mathbb P^1$ embedded with $2B$ is a minimal degree variety; see \cite[Example 5.6]{BSV}. Our approach also applies to this situation whenever $d$ or $n$ is odd, since the toric catalecticant $\phi_{f,B}$ of a general $f\in S_{2B}$ has full rank $(d+1)(n+1)$ and this is even; the last restriction is required in the proof of Proposition \ref{changeofsignAntip} so that the total degree of $\Omega(f)(\ell)$ is odd in the coefficients of $f$. 
\end{remark}

\section{Real rank boundary and hyperdeterminants}

In this section, we continue the discussion of real rank boundaries and  connect them to hyperdeterminants; we refer to \cite{GKZ} for the theory of hyperdeterminants. In the following, we identify a projective space with its dual. 

\begin{proposition}\label{hypervstang}
Let $X =\mathbb P^1\times \mathbb P^{n-1} \times \mathbb P^{n-1}$. Then the join of its tangential variety and its $(n-2)$nd secant variety, $\tau(X)+\sigma_{n-2}(X)$, coincides with the hyperdeterminant, the dual variety $X^{\vee}$. 
\begin{proof}

Let us consider the space $\mathbb P(\mathbb C^2\otimes \mathbb C^n \otimes \mathbb C^n)$ as a space of pencils of $n\times n$ matrices. 
General tensors in $\mathbb P(\mathbb C^2\otimes \mathbb C^n \otimes \mathbb C^n)$ have the form $e_1\otimes \textnormal{Id}_n + e_2\otimes M$, where $M$ is a diagonal matrix with distinct generic eigenvalues, as a general pencil contains an invertible matrix and generic matrices are diagonalizable. (Generic matrices are diagonalizable because generic univariate polynomials have distinct complex roots.) The general point in $\tau(X)+\sigma_{n-2}(X)$ is of the form 
\begin{equation}\label{tensorintang}
a_1\otimes \textnormal{Id}_n+a_2\otimes \begin{pmatrix} 
\lambda_1 & 1 & 0 & \ldots&0 \\
0 & \lambda_1 &0 & \ldots& 0 \\
0& 0 & \lambda_2 & \ldots & 0 \\
\ldots & \ldots &\ldots &\ldots &\ldots \\
0 & 0 & 0 & \ldots & \lambda_{n-1} \\
\end{pmatrix}, 
\end{equation}
\noindent where $\textnormal{Id}_n$ is the identity $n\times n$ matrix. Indeed, this point can be written as 
$$
(a_1+\lambda_1 a_2)\otimes a_1\otimes a_1+(a_1+\lambda_1 a_2)\otimes a_2\otimes a_2 +a_2\otimes a_1\otimes a_2 + \sum_{i=2}^{n-1} (a_1+\lambda_i a_2)\otimes a_{i+1}\otimes a_{i+1},
$$ 
\noindent where the first summand is in the tangent space of $X$ at $(a_1+\lambda_1a_2)\otimes a_1 \otimes a_2$. Hence this point lies in $\tau(X)+\sigma_{n-2}(X)$. Moreover, by Schl\"{a}fli's method \cite[Chapter 14]{GKZ}, and since the corresponding binary form of $(\ref{tensorintang})$ has a double root, this point is in $X^{\vee}$. Since both of $\tau(X)+\sigma_{n-2}(X)$ and $X^{\vee}$ are irreducible hypersurfaces, the statement follows.
\end{proof}
\end{proposition}

\begin{remark}
Proposition \ref{hypervstang} holds for tensors with symmetric $n\times n$ slices. These can be seen as points in $\mathbb P(\mathbb C^2\otimes S^2 \mathbb C^n)$. The hyperdeterminant in this space is a linear section of the hyperdeterminant in $\mathbb P(\mathbb C^2\otimes \mathbb C^n \otimes \mathbb C^n)$. In this case, the Jordan block $\begin{pmatrix} \lambda_1 & 1 \\ 
0  & \lambda_1 \\ \end{pmatrix}$ above is replaced by the symmetric block $\begin{pmatrix} 1+\lambda_ 1& i \\ 
i  & -1+\lambda_1 \\ \end{pmatrix}$. The consequence is that the hyperdeterminant in $\mathbb P(\mathbb C^2\otimes S^2 \mathbb C^n)$ is the join of $\tau(X)$ and $\sigma_{n-2}(X)$, where $X$ is the Segre variety $\mathbb P^1\times \mathbb P^{n-1}$ embedded with the line bundle $B = (1,2)$. 
\end{remark}

\begin{remark}
The symmetric analogue of Proposition \ref{hypervstang} for $n=2$ is well-known. The tangential variety of the twisted cubic $\nu_3(\mathbb P^1)\subset \mathbb P^3 = \lbrace \mathcal{C} = a_0x^3+a_1x^2y+a_2xy^2+a_3y^3 \rbrace$ is a quartic surface (classically known as {\it tangent developable}) given by cubics with a double root. It coincides with the zero set of the discriminant 
$$
\textnormal{Disc}(\mathcal{C}) = a_1^2a_2^2 - 4a_0a_2^3 - 4a_1^3a_3 - 27a_0^2a_3^2 + 18a_0a_1a_2a_3.
$$
\noindent The discriminant is the specialization of the hyperdeterminant to symmetric tensors. 
\end{remark}

For $X = \mathbb P^1\times \mathbb P^1 \times \mathbb P^1$, as shown by De Silva and Lim \cite[Table 7.1]{DeSLim}, the hyperdeterminant constitutes the real rank boundary between the only two possible typical ranks, namely two and three.  More generally, for $X = \mathbb P^1\times \mathbb P^{n-1} \times \mathbb P^{n-1}$, there are only two typical ranks, $n$ and $n+1$. In  \cite[Theorem 1]{GBerq}, Bergqvist proved the following result. 

\begin{theorem}[{\bf  \cite[Theorem 1]{GBerq}}]\label{Bergqvist}
Let $n\geq 2$. Let $T$ be a general tensor in $\mathbb P(\mathbb R^2\otimes \mathbb R^n \otimes \mathbb R^n)$. Let $T_1$ and $T_2$ be the two slices of $T$ with respect to the basis $a_1,a_2$ of $\mathbb R^2$. Then $T$ has real rank $n$ if and only if the binary form $p_{T}(a_1,a_2)=\det (a_1T_1+a_2T_2)$ has $n$ real solutions. 
\end{theorem}

The following generalizes the result of De Silva and Lim.

\begin{corollary}
The real rank boundary between the typical ranks $n$ and $n+1$ is the hyperdeterminant. 
\begin{proof}
By Schl\"{a}fli's method \cite[Chapter 14]{GKZ}, the hyperdeterminant of $X= \mathbb P^1\times \mathbb P^{n-1} \times \mathbb P^{n-1}$ is the complexification of the discriminant of the binary form $p_{T}(a_1,a_2)$, where $T$ is a tensor in $\mathbb P(\mathbb R^2\otimes \mathbb R^n \otimes \mathbb R^n)$. In each connected component of the complement of the discriminant of $p_{T}(a_1,a_2)$, the number of real solutions is constant. Hence the number of real solutions of $p_{T}(a_1,a_2)$ changes whenever we cross the discriminant. By Theorem \ref{Bergqvist}, a general real tensor $T$ has real rank $n$ if and only if the corresponding $p_{T}(a_1,a_2)$ has $n$ real solutions. Suppose that these $n$ real solutions are all distinct. Then $T$ is in a connected component $\mathcal C_0$ of the complement of the discriminant. As we move $T$ and cross the discriminant, the number of real roots of $p_{T}(a_1,a_2)$ changes. Hence we reach another connected component $\mathcal C_1$, where the number of {\it real} roots of $p_{T}(a_1,a_2)$ is strictly less than $n$. By Theorem \ref{Bergqvist}, in $\mathcal C_1$, there exist tensors whose real rank is greater than or equal to $n+1$. Since $n+1$ is the only typical rank larger than $n$, in $\mathcal C_1$ there exists an open set of tensors of real rank $n+1$. Thus the hyperdeterminant coincides with the real rank boundary. 
\end{proof}

\end{corollary}

As a corollary of Proposition \ref{hypervstang}, we derive: 

\begin{corollary}\label{newinterpofbound}
The real rank boundary of $X = \mathbb P^1\times \mathbb P^{n-1} \times \mathbb P^{n-1}$ is the join of its tangential and its $(n-2)$nd secant variety.
\end{corollary}

\begin{remark}
Proposition \ref{hypervstang} provides a class of examples supporting \cite[Conjecture 5.5]{MMSV16}, because in such cases we have generic identifiability.  
\end{remark}

\begin{remark}
By \cite[Proposition 2.7]{MMSV16}, the hyperdeterminant of $X = \mathbb P^1\times \mathbb P^{n-1} \times \mathbb P^{n-1}$ is the Hurwitz form of 
the determinantal variety of $n\times n$ matrices of rank at most $n-1$. This gives another interpretation of the real rank boundary of $X$. 
\end{remark}

One may wonder whether in all the cases, where Schl\"{a}fli's method applies, the hyperdeterminant has the peculiar interpretation as a join of the tangential and an appropriate secant variety. This is not the case as shown by the next two examples.  

\begin{example}[{\bf The $3\times 3\times 3$ hyperdeterminant}]
Let $X=\mathbb P^2\times \mathbb P^2\times \mathbb P^2$. The variety $X$ is defective, as its fourth secant is a hypersurface of degree nine. Its fifth secant $\sigma_5(X)$ fills the ambient space $\mathbb P^{26}$. Consider the join $Y=\tau(X)+\sigma_3(X)$. This does not coincide with the $3\times 3\times 3$ hyperdeterminant. Indeed, let us consider the point $T = (\lambda_{0,1}a_0+\lambda_{1,1}a_1)\otimes b_2\otimes c_2 + a_2\otimes (\lambda_{0,2}b_0+\lambda_{1,2}b_1)\otimes c_2+a_2\otimes b_2\otimes(\lambda_{0,3}c_0+\lambda_{1,3}c_1)+a_0\otimes b_0\otimes c_0+a_1\otimes b_1\otimes c_1+a_2\otimes b_2\otimes c_2$. The first summand is a point in the tangent space to $X$ at $a_2\otimes b_2\otimes c_2$ and the second is a point in the third secant to $X$. Using Schl\"{a}fli's method, the corresponding polynomial in the variables $a_i$ is $p_T(a_0,a_1,a_2) = \lambda_{0,1}a_0^2a_1+\lambda_{1,1}a_0a_1^2-\lambda_{1,2}\lambda_{1,3}a_0a_2^2-\lambda_{0,2}\lambda_{0,3}a_1a_2^2+a_0a_1a_2$. This is a smooth ternary cubic for a generic choice of the coefficients. Thus the point $T$ does not lie in the $3\times 3\times 3$ hyperdeterminant. 
\end{example}

\begin{example}[{\bf The $2\times 2\times 2\times 2$ hyperdeterminant}]
The hyperdeterminant of $X=\mathbb P^1 \times \mathbb P^1 \times \mathbb P^1 \times \mathbb P^1$ is computed analogously as above, namely using a version of Schl\"{a}fli's method. We briefly recall this here. Let us denote the coordinates of $\mathbb P(\mathbb C^2\otimes \mathbb C^2 \otimes \mathbb C^2 \otimes \mathbb C^2)$ by 
$z_{i,j,k,l}$, for $i,j,k,l\in \lbrace 0,1\rbrace$. Consider the $2\times 2 \times 2$ hyperdeterminant as a polynomial in coordinates $u_{i,j,k}$. We substitute $u_{i,j,k}$ with
$z_{i,j,k,0}+z_{i,j,k,1}w$. We obtain a degree four polynomial $p(w)$ with coefficients in $z_{i,j,k,l}$. The discriminant of $p(w)$ is the $2\times 2 \times 2\times 2$ hyperdeterminant. We do not have the equality between the hyperdeterminant of $X$ and the join variety $Y=\tau(X)+\sigma_2(X)$. Indeed, let us pick the point $T = a_1\otimes b_0 \otimes c_0\otimes d_0 + a_0\otimes b_1\otimes c_0\otimes d_0+a_0\otimes b_0 \otimes c_1\otimes d_0 + a_0\otimes b_0\otimes c_0\otimes d_1+a_0\otimes b_0\otimes c_0\otimes d_0 + a_1\otimes b_1\otimes c_1 \otimes d_1$ in $Y$. The corresponding polynomial is $p_T(w) = w^4+2w^3+w^2+4w$. Since this polynomial has not double roots, its discriminant is not zero. Thus the point $T$ does not lie in the $2\times 2\times 2\times 2$ hyperdeterminant. 
\end{example}

It is natural to ask the ensuing

\begin{question}
Are there other instances where the hyperdeterminant is a join between the tangential and the appropriate secant variety? What about symmetric cases?
\end{question}

\medskip

\begin{small}

\noindent
{\bf Acknowledgements.}\smallskip \\
\noindent I  would like to thank Mateusz Micha\l{}ek and Bernd Sturmfels for their constant support and for inspiring discussions. I would like to thank Mateusz Micha\l{}ek for insightful comments to an earlier version of this article. I would like to thank Alessandra Bernardi, Alessandro Oneto, Giorgio Ottaviani, Anna Seigal for very useful conversations on topics related to this work, and Zach Teitler for sending me a copy of Kleppe's thesis. I would like to thank the Max Planck Institute for Mathematics in the Sciences in Leipzig for providing a wonderful working environment. Finally, I would like to thank the referee for very useful remarks and suggestions. 

\end{small}
\medskip

\begin{small}

\end{small}

\bigskip

\noindent
\footnotesize {\bf Author's address:}\\
\noindent 
Texas A\&M University, \\
Department of Mathematics\\
Mailstop 3368\\
College Station, TX 77843-3368\\
{\tt eventura@math.tamu.edu}
\end{document}